\documentstyle[12pt]{article}
\begin{document}

\textwidth18cm
\textheight20cm

\newtheorem{Ex}{\bf Example}
\newtheorem{Rem}{\bf Remark}
\newtheorem{Sa}{\bf Proposition}
\newtheorem{Le}[Sa]{\bf Lemma}
\newtheorem{Hi}[Sa]{\bf Lemma}
\newtheorem{Th}[Sa]{\bf Theorem}
\newtheorem{Fo}[Sa]{\bf  Corollary}
\newtheorem{Beh}[Sa]{\bf Proposition}
\newtheorem{De}{\bf Definition}

\newcommand{\DS}{ \mathrel{\mathop{+}\limits^.   }}
\newcommand{\RR}{{\rm I\kern-0.14em R }}
\newcommand{\PP}{{\rm I\kern-0.14em P }}
\newcommand{\NN}{{\rm l\kern-0.14em N }}
\newcommand{\CC}{{\rm\raise 0.001ex\vbox{\hrule height
1.22ex width 0.8pt}\kern-0.29em C}}
\newcommand{\ZZ}{ {\sf Z}\hspace{-0.4em}{\sf Z}\ }
\newcommand{\TT}{ { /}\hspace{-0.4em}{|} \ }
\newcommand{\AAAAA}{ \hspace{-0.6em}{\bf a}\ }
\def\ii{{\rm i}}
\def\dd{{\, \rm d }}


\begin{center}
{\LARGE   Quantum theory of the real and the complexified
projective line }
\end{center}

\begin{flushleft}
{\large FRANK LEITENBERGER } \\
{\small \it Fachbereich Mathematik, Universit\"at Rostock,
Rostock, D-18051, Germany. \\
e-mail: frank.leitenberger@mathematik.uni-rostock.de }
\end{flushleft}


\vspace{0.5cm}
{\small
{\bf Abstract: \rm   Quantum deformations
of sets of points of the  real and the
complexified projective line are constructed.
These deformations depend on the deformation
parameter $q$ and certain further parameters $\lambda_{ij}$.
The deformations for which
the subspace of polynomials of degree three
has a basis of ordered monomials are selected.
We show that the corresponding algebras of three
points have "polynomiality".
Invariant elements which
turn out to be cross ratios in the classical
limit  are defined. For the  special case
$|\lambda_{ij}| =  1 $ a quantum cross ratio
with properties similar to the classical case is presented.
As an application a quantum version of
the real Euclidean distance is given.}}

\vspace{0.5cm}


\newpage
\noindent
{\Large \bf I. Introduction }
\\ \\
Noncommutative geometry can be viewed as the quantum deformation
of an increasing set of mathematical objects and theories.
In this article we introduce a quantized version of the
projective geometry in one dimension.

First we consider the algebra of coordinates of points
$t_i$, $i\in I$ of the real line.
All possible Poisson brackets between these variables
(which are covariant with respect
to the action of the standard
Poisson structure of the M\"obius group) are parametrisized
by a large set of parameters $\lambda_{ij}$, $i,j\in I$
which obey certain second order polynomial conditions.
These quasiclassical structures give a good
insight into the quantized situation.

In section 3 we replace the algebra of coordinates $t_i$
by a $U_q ( sl(2,\RR )$-modul algebra of noncommutative
coordinates $v_i$. We show that (up to reordering the points)
third degree polynomials have a basis of ordered polynomials
if a condition in analogy to the Jacobi identity is satisfied.
There is one exceptional structure, which has no
quasiclassical limit.

For algebras of three points we give a proof
that the subspaces of homogeneous polynomials of higher degrees
have also the classical dimension.
The condition (which plays the role of the
Jacobi identity) takes the form of the Einstein addition law
for velocities, where the light velocity corresponds to
the quantization parameter $q$.
We consider also the complexified case.

In section 4 we consider invariant elements
$C_{ijkl}$ of four points,
which turn out to be cross ratios in the classical limit.
Invariant elements occured first in Ref. 1 in a similar
situation for the quantum sphere.
Our class of invariant elements contain the elements
of the type of Ref. 1 as a subclass.
However these elements do not satisfy a reality condition
($*$-invariance) and have no simple behaviour with respect
to permutations of the four points.
We show that for parameters $|\lambda_{ij}|= 1$ certain
rational functions $f(C_{ijkl},q)$ have the classical
properties of the cross ratio.

As an application we define a quantization of the Euclidean
distance on the real line.

\newpage
\noindent
{\Large \bf II. The quasiclassical limit }
\\ \\
\noindent
Consider the one-dimensional real projective line $\PP$, the action
of the M\"obius group $G$ on it and an arbitrary set of variable
points $P_i\in \PP, i\in I$ on it.
Let $t$ be a coordinate on $\PP$ and $t_i$ the coordinates of the
points $P_i$ and let
the M\"obius group
$G\cong SL(2,\RR)/ \{\pm I \}$
act on the points  $P_i$ according to
\[   t_i \longrightarrow  \frac{at_i+b}{ct_i+d}.    \]
Endow $G$ with the structure of a Poisson group:
\[   \begin{array}{ccc}
 \{a,b\}  =   ab,  \ \ \ &
 \{a,c\}  =   ac,  \ \ \ &
 \{a,d\}  =   2 bc, \\
 \{b,c\}  =    0, \ \ \  &
 \{b,d\}  =   bd,  \ \ \ &
 \{c,d\}  =   cd.
\end{array} \]
Let $I$ be an arbitrary set.
The general $G$-covariant (cf. Ref. 2) pre-Poisson bracket
on the algebra of polynomials
of the variables $t_i,\ \ i\in I$ is given by
\begin{eqnarray}
\{ t_j, t_i \} = t_j^2 - t_i^2 - \lambda_{ij}
                 (t_i - t_j)^2
\end{eqnarray}
with arbitrary coefficients $\lambda_{ij}=-\lambda_{ji}$.

\noindent
{\it Proof.} The G-covariant pre-Poisson brackets are given
as sums of a special
G-covariant bracket and an arbitrary G-invariant bracket
(cf. Ref. 2). The general G-invariant bracket is given by
\[ \{ t_j, t_i \} =  \lambda_{ij}  (t_i - t_j)^2.   \]
Therefore it is sufficient to check the G-covariance for
one of the above structures by explicite calculation .
$\bullet$

A calculation yields that
this pre-Poisson bracket becomes a Poisson bracket
(i.e. it satisfies the Jacobi identity) if and only if
\begin{eqnarray}
   \lambda_{ij}\lambda_{jk}-\lambda_{ij}\lambda_{ik}
   - \lambda_{ik}\lambda_{jk} = -1 \end{eqnarray}
for arbitrary three different $i,j,k \in I$. The relation is
invariant with respect to a permutation of the indices.

Instead to give the general solution of this equation
we prefer to describe some basic examples.

\noindent
{\it Example $1$.} Let $I$ be a linear ordered set.
$I$ may have arbitrary cardinality.
Set $\lambda_{ij}=1$ for $i<j$, $i,j\in I$.
We obtain
\[ \{ t_j, t_i \} = 2 t_i(t_j - t_i), \ \ \ (i<j). \]
Conversely, let $|\lambda_{ij}|= 1$ for arbitrary different
$i,j\in I$.
Then we can define a linear order relation on $I$.
We set $i<j$, if $\lambda_{ij}=1$
(the transitivity law of the order relation
corresponds to the Jacobi identity).

\noindent
{\it Example $2$.} Let $I=\ZZ$ be the set of integers,
$( \alpha_k )_{ k\in \ZZ }$ be a sequence of positive
integers and let
\[ \lambda_{ij}= \coth ( \sum_{k=i}^{j-1} \alpha_k ),\ \ \ \ \ \
\alpha_k > 0,\ \ \ i<j. \]
The Poisson property follows from the addition formula
\[ \coth (\alpha +\beta ) =
   \frac{ 1 + \coth( \alpha ) \coth ( \beta ) }{
   \coth( \alpha ) + \coth( \beta )}. \]
We have
\[ \{ t_j, t_i \} = t_j^2 - t_i^2 - \coth ( \sum_{k=i}^{j-1} \alpha_k )
  (t_i - t_j)^2,\ \ \ i<j. \]
Note that $|\lambda_{ij}|> 1$ because of $|\coth(\alpha)|>1$ for
$\alpha\neq 0$.
We obtain again Example 1 in the limit $\alpha_k \rightarrow \infty$.

\noindent
{\it Example $3$.}
In order to get examples with $|\lambda_{ij}| < 1$
for certain $i,j\in I$, let $I=\ZZ$ and
\[ \lambda_{ij}=
\coth ( \sum_{k=i}^{j-1} (\alpha_k + \frac{i\pi }{2}))
,\ \ \ \ \ \
\alpha_k > 0. \]
Then we have $|\lambda_{ij}| < 1$ for $i-j$ odd and
        $|\lambda_{ij}| > 1$ for $i-j$ even.

\noindent
{\it Example $4$.} Now we give an example with
$\lambda_{ij}=0$ for certain $i,j\in I$.
We remark that
if $i$ is fixed, there is only one different index
$j$ with $\lambda_{ij}=0$ (because of equation (2)).
Let $I=\ZZ \times {\ZZ}_2$ and
$z_i:=v_{(i,0)}$, $z_{ {\overline{i}} }^*:= v_{(i,1)}$
be the corresponding variables. One can check equation
(2) for $\lambda_{i,{\overline{i}}}=0$ and
$\lambda_{ij}=
 \lambda_{{\overline{i}}j}=
 \lambda_{i{\overline{j}}}=
 \lambda_{{\overline{i}}{\overline{j}}}=1$
for $i<j$.
The Poisson bracket is given by
\[ \{ z_{\overline{i}}^*,z_{i} \}=
                     (z_{\overline{i}}^{*2} - z_{i}^{2}),   \]
\[ \{ z_j  ,z_i     \}=
           2 z_i (z_j-z_i),\ \ \ i<j,            \]
\[ \{ z_{\overline{j}}^*,
      z_{\overline{i}}^*   \}=
           2 z_{\overline{i}}^*
           (z_{\overline{j}}^{*}-z_{\overline{i}}^*),\ \ \ i<j,    \]
\[ \{ z_j  ,z_{\overline{i}}^*   \}=
           2 z_{\overline{i}}^*
           (z_{j} - z_{\overline{i}}^*),\ \ \ i<j,      \]
\[ \{ z_{\overline{j}}^*,z_i     \}=
           2 z_i (z_{\overline{j}}^*  - z_i ),\ \ \ i<j.       \]
Because $z_i$ and $z_i^*$ transform according to
\[   z_i   \longrightarrow  \frac{az_i  +b}{cz_i  +d},  \ \ \ \ \
     z_i^* \longrightarrow  \frac{az_i^*+b}{cz_i^*+d}          \]
we identify $z_i$ and $z_i^*$ with complex and conjugate
complex coordinates of the complexified projective line.

\noindent
{\it Example $5$.} Now  we explain, how to put structures
together to a new one. Let $I'$ be an ordered set.
Consider variables $t_i^{(j)}$, $i\in I^{(j)}$, $j\in I'$
(i.e. $I=\bigcup_{j\in I'} I^{(j)}$).
Suppose that for fixed $j$ the variabels
$t_i^{(j)}$, $i\in I^{(j)}$ carry
an arbitrary Poisson structure. We combine these structures
according to
 \[ \{ t_{i'}^{(j')}  , t_{i}^{(j)} \}
   = 2  t_{i}^{(j)} ( t_{i'}^{(j')}  - t_{i}^{(j)}  ),
   \ \ \ j<j'. \]
The verification
uses the fact that equation (2) is
satisfied, if two of the numbers
$\lambda_{ij}$, $\lambda_{ik}$
and $\lambda_{jk}$
are equal to 1.

\noindent
{\it Example $6$.}  An example with complex
$\lambda_{ij}$ is given by  $I=\{1,2,3,4 \}$ and
$\lambda_{12}=\lambda_{23}=\lambda_{34}={\rm i}$,
$\lambda_{14}=- {\rm i}$
and $\lambda_{13}=\lambda_{24}=0$.

\newpage
\noindent
{\Large \bf III. Quantization of point sets of the real
and complex projective line }
\\ \\
\noindent
{\bf A. A class of $U_q (sl_2,\RR )$-modul algebras} \\
Let $U_q (sl_2,\RR ),\ |q|=1$, $q\neq \pm 1,\pm i$
be the unital $*$-algebra, given by
generators $E, F, K, K^{-1}$ and relations
\[    KK^{-1}=K^{-1}K=I,
 \ \  KE=qEK,
 \ \  KF=q^{-1}FK                              \]
\[     [E,F]=\frac{ K^2-K^{-2} }{ q-q^{-1} },  \]
and the involution
\[ E^* = F,\ \ \ \ \ F^* = E,\ \ \ \ \  K^* = K^{-1}. \]

For an arbitrary set $I$
consider the unital
complex algebra $A$
which is freely generated
by the elements $x_i,\ i\in I$.

An action of $U_q(sl(2,\RR))$ on $A$ is given,
if we set
\begin{eqnarray}
E(1)=F(1)=0,\ K(1)=1,\ \ E(x_i)=1,\ F(x_i)=-x_i^2,\
                        K(x_i)=q^{-1}x_i
\end{eqnarray}
and require for $x,y \in A$
\begin{eqnarray} K(xy) = K(x) K(y)              \end{eqnarray}
\begin{eqnarray} E(xy) = E(x)K(y)+K^{-1}(x)E(y),\end{eqnarray}
\begin{eqnarray} F(xy) = F(x)K(y)+K^{-1}(x)F(y).\end{eqnarray}
\noindent

A proof can be given along the lines of Ref. 3,
where a similar modul of two variables occurs.

Further we consider ideals $I$ of $A$, which are generated by
all (non\-vani\-shing) elements of the form
\begin{eqnarray}  X_{ij}:=
  a_{ij} x_i   x_j + b_{ij} x_i^2
+ c_{ij} x_j^2     + d_{ij} x_j   x_i,
\end{eqnarray}
where $a_{ij}=d_{ji}$, $b_{ij}=c_{ji}$
(i.e. the elements are invariant with respect to a change
of $i$ and $j$).
The action of $U_q(sl(2,\RR))$ on $A$ induces an action
on the factor algebra $A/I$ if and only if $I$ is
$U_q(sl(2,\RR))$-invariant. We denote the image of the
quotient map of $x_i$ by $v_i$.

\begin{Sa}  The ideal $I$ is $U_q (sl(2, \RR))$-invariant
if and only if
\[    a_{ij} q^2 + b_{ij} (1+q^2) + d_{ij} = 0 \]
and
\begin{eqnarray}
       a_{ij}  +  c_{ij} (1+q^2) + d_{ij} q^2 = 0.
  \end{eqnarray}
The ideal $I$ is generated by the elements
\[ I_{ij}^\lambda = [x_j,x_i]-\frac{q^2-1}{q^2+1}
             (x_j^2-x_i^2 - \lambda_{ij}(x_i-x_j)^2), \]
with $\lambda_{ij}=
 \frac{q^2+1}{q^2-1}
 \frac{ a_{ij}+d_{ij} }{ a_{ij} - d_{ij} },$
($\lambda_{ij}=-\lambda_{ji}$, $i,j\in I$)
for $a_{ij} \neq d_{ij}$
and
\[ I_{ij}^\infty =  (x_i-x_j)^2. \]
for $a_{ij} = d_{ij}$.
\end{Sa}

\noindent {\it Proof:}
First we consider the case $a_{ij}\neq d_{ij}$.
Using (3), (4), (5), (6), we have
\[ E(X_{ij}):= E(a_{ij} x_i   x_j + b_{ij} x_i^2
   + c_{ij} x_j^2     + d_{ij} x_j   x_i)\]
\[ =   ( a_{ij} q  +  b_{ij} (q+\frac{1}{q}) +  d_{ij} \frac{1}{q} ) x_i
 + ( a_{ij} \frac{1}{q}  +  c_{ij} (q+\frac{1}{q}) + d_{ij} q ) x_j. \]
Because there are no linear relations between the generators,
the last expression is contained in $I$, if and only if it is zero.
Therefore the equations (8) are necessary.

Now we show, that the equations (8) are also sufficient.
We have to show, that $K(X_{ij}),F(X_{ij}) \in I$.
We have
\[  K(X_{ij}):= q^{-2} X_{ij} \in I.  \]
Further
\[  F(X_{ij}):= -  \frac{a_{ij}}{q} x_i^2 x_j - a_{ij} q x_i x_j^2
                -  \frac{b_{ij}}{q} x_i^3     - b_{ij} q x_i^3
                -  \frac{c_{ij}}{q} x_j^3     - c_{ij} q x_j^3
                -  \frac{d_{ij}}{q} x_j^2 x_i - d_{ij} q x_i x_j^2.  \]
Now let $\alpha :=\frac{1+q^2}{q} \frac{a_{ij}}{d_{ij}-a_{ij}} $,
        $\beta  :=\frac{1+q^2}{q} \frac{d_{ij}}{a_{ij}-d_{ij}} $.
Because of (8) we have
\[ \alpha(a_{ij}+b_{ij}) = -\frac{a_{ij}}{q},\ \
   \alpha(a_{ij}+c_{ij}) = - a_{ij} q,\ \
   \beta (c_{ij}+d_{ij}) = -\frac{d_{ij}}{q},\ \
   \beta (b_{ij}+d_{ij}) = - d_{ij} q, \]
\[ (\alpha+\beta)  = - (q+\frac{1}{q}),\ \
   (d_{ij}\alpha+a_{ij}\beta)  = 0.       \]
Therefore
\[ F(X_{ij})=   \alpha (a_{ij}+b_{ij}) x_i^2 x_j
              + \alpha (a_{ij}+c_{ij})x_i x_j^2
              + \beta  (c_{ij}+d_{ij})x_j^2 x_i
              + \beta  (b_{ij}+d_{ij}) x_j x_i^2 + \]
\[            + (\alpha+\beta )b_{ij} x_i^3
              + (\alpha+\beta )c_{ij} x_j^3  +
      (d_{ij}\alpha+a_{ij}\beta)x_i x_j x_i
     +(d_{ij}\alpha+a_{ij}\beta)x_j x_i x_j \]
\[ = \alpha (x_i X_{ij}+X_{ij} x_j)+ \beta (x_j X_{ij}  +X_{ij} x_i)
\in I. \]
I.e. the equations (8) are sufficient.

Further we confirm the second statement
of the proposition. We show that the generators
$X_{ij}$ and $I_{ij}^\lambda$ differ only by a
constant nonzero factor.
{\small
\[ I_{ij}^\lambda =
  ( -1 - \lambda_{ij} \frac{q^2-1}{q^2+1} ) x_i x_j
+ ( 1+\lambda_{ij}) \frac{q^2-1}{q^2+1}      x_i^2
+ (-1+\lambda_{ij}) \frac{q^2-1}{q^2+1}      x_j^2
+ (  1 - \lambda_{ij} \frac{q^2-1}{q^2+1} ) x_j x_i        \]
\[ = \frac{ -2 }{ a_{ij}-d_{ij} }
 ( a_{ij} x_i x_j
 + \frac{ -a_{ij}q^2-d_{ij} }{1+q^2} x_i^2
 + \frac{ -a_{ij}-d_{ij}q^2 }{1+q^2} x_j^2
 + d_{ij} x_j x_i ) \]
\[ = \frac{-2 }{ a_{ij}-d_{ij} }
     (  a_{ij} x_i x_j
      + b_{ij} x_i^2  + c_{ij} x_j^2 + d_{ij} x_j x_i  ).  \]}
The second statement of the proposition follows.

Finally, let $a_{ij}=d_{ij}\neq 0$. We obtain
$b_{ij}=c_{ij}=-a_{ij}=-d_{ij}$. The proposition follows.
$\bullet$

We exclude the case $I_{ij}^\infty$
from the following considerations.
Proposition 1 admits the following Definition.

\begin{De} \rm Let $I$ be an arbitrary set.
${A'}_\lambda$, $\lambda=(\lambda_{ij};i,j\in I)$
 is the left $U_q(sl_2,\RR )$-modul algebra given by generators
 $v_i,\ i\in I$ and relations
\begin{eqnarray}
[v_j,v_i]=\frac{q^2 -1}{q^2+1}
             (v_j^2-v_i^2-\lambda_{ij}(v_i-v_j)^2)
\end{eqnarray}
($\lambda_{ij}=-\lambda_{ji}$, $i,j\in I$).
\end{De}
Let $\lambda_{ij}$ be fixed.
In the quasiclassical limit $q\rightarrow 1$ the bracket (9)
becomes the pre-Poisson bracket (1).
By Proposition 1 the moduls ${A'}_\lambda$
are the most general
$U_q(sl(2,\RR))$-moduls given by generators
and quadratic relations.

In order to define crossratios (see below) we need also
elements of the form $(\sum \alpha_i v_i)^{-1}$.

\begin{De}\rm We obtain the extension $A_\lambda$ of
${A'}_\lambda$ if we adjoin the inverses of
the nonvanishing finite sums
$(\sum \alpha_i v_i)^{-1}$  to $A_\lambda$ and
require that
\begin{eqnarray}
(\sum \alpha_i v_i)^{-1}(\sum \alpha_i v_i)
 =(\sum \alpha_i v_i)(\sum \alpha_i v_i)^{-1}=1.
\end{eqnarray}
The $U_q(sl_2,\RR)$-modul structure extends by the
requirements (2), (3), (4) to all $x,y\in A_\lambda$.
\end{De}

\vspace{0.5cm}
\noindent
{\bf B. Polynomials of degree three}

\noindent
There is a principle about
a correspondence between the Jacobi identity
for Poisson brackets and the Poincar\'{e}-Birkhoff-Witt-Theorem for
the quantized algebras.
We will prove the statement of the PBW-Theorem only for polynomials
of degree three.
We obtain that condition (2) for the Jacobi identity characterizes
also the classical behaviour of third order polynomials
(with the exception of two special values for
$\lambda_{i,j}$ and one structure where (2) is not satisfied).
The main part of the proof is Lemma 3.

\noindent
{\bf Remark:}
Let $b_{ij},c_{ij}\neq 0$ (i.e. $|\lambda_{ij}|\neq 1$).
Then on the right side in relation (12)
occurs a  monomial which is lexicographically smaller and
a monomial which is lexicographically greater than the monomial
on the left side. Therefore we can not use the Bergman theory
(cf. Ref. 4) for a proof of the PBW-Theorem.



Let $i<j$ and suppose $d_{ij}\neq 0$
(or the equivalent condition
$\lambda_{i,j} \neq  \frac{q^2+1}{q^2-1}$).
It follows
\begin{eqnarray}
 v_j v_i := \alpha_{ij} v_i v_j
          + \beta_{ij}  v_i^2
          + \gamma_{ij} v_j^2
\end{eqnarray}
with
\[\alpha_{ij}:=\frac{-a_{ij} }{d_{ij}}=
\frac{ -(q^2+1)-\lambda_{ij}(q^2-1) }{\lambda_{ij}(q^2-1)-(q^2+1)},\]

\[\beta_{ij}:=\frac{-b_{ij} }{d_{ij}}
                  =\frac{1-q^2 \alpha_{ij}}{1+q^2}=
\frac{(q^2-1)(1+\lambda_{ij})}{\lambda_{ij}(q^2-1)-(q^2+1)},\]

\[\gamma_{ij}:=\frac{-c_{ij} }{d_{ij}}
                  =\frac{q^2 - \alpha_{ij}}{1+q^2}=
\frac{(q^2-1)(-1+\lambda_{ij})}{\lambda_{ij}(q^2-1)-(q^2+1)}.\]
(cf. formulas (7),(8),(9)).

\begin{Le}
(i) Let $i<j$, $\lambda_{i,j} \neq \frac{q^2+1}{q^2-1}$
and $v_i,v_j$ two of the generators of ${A'}_\lambda$.
Then
\[ (1-\beta_{ij} \gamma_{ij})                           v_j^2 v_i =
   ( \beta_{ij}^2            ( 1 + \alpha_{ij} )      ) v_i^3
 + ( \alpha_{ij} \beta_{ij}  ( 1 + \alpha_{ij} )      ) v_i^2 v_j
 + ( \alpha_{ij}^2 + \alpha_{ij}\beta_{ij}\gamma_{ij} ) v_i v_j^2
 + ( \gamma_{ij}             ( 1 + \alpha_{ij} )      ) v_j^3.
\]
We have $1-\beta_{ij} \gamma_{ij} \neq 0$ for the coefficient
on the left side if and only if
$\lambda_{i,j} \neq  \frac{q^4+1}{q^4-1}$.

(ii) Let $1-\beta_{ij} \gamma_{ij} \neq 0$. Then we can reduce
every monomial of degree three into a sum of ordered monomials.
\end{Le}

\noindent
{\it Proof.}
(i) We have
\[  v_j v_j v_i =   \alpha_{ij} v_j v_i v_j
                   + \beta_{ij}  v_j v_i^2
                   + \gamma_{ij} v_j^3       \]
\[              =   \alpha_{ij} v_j v_i v_j +
    (\beta_{ij}  \alpha_{ij} v_i v_j v_i    +
     \beta_{ij}  \beta_{ij}  v_i^3          +
     \gamma_{ij} \beta_{ij}  v_j^2 v_i)     +
                 \gamma_{ij} v_j^3           \]
\[ =(\alpha_{ij} \alpha_{ij} v_i v_j^2   +
     \alpha_{ij} \beta_{ij}  v_i^2 v_j   +
     \alpha_{ij} \gamma_{ij} v_j^3)      +
    (\beta_{ij}  \alpha_{ij} \alpha_{ij} v_i^2 v_j  + \]
\[   \beta_{ij}  \alpha_{ij} \beta_{ij}  v_i^3      +
     \beta_{ij}  \alpha_{ij} \gamma_{ij} v_i v_j^2) +
     \beta_{ij}  \beta_{ij}  v_i^3       +
     \gamma_{ij} \beta_{ij}  v_j^2 v_i   +
                 \gamma_{ij} v_j^3.           \]
The identity follows. Further
\[ (1-\beta_{ij} \gamma_{ij})
   =1-\frac{(q^2-1)( 1+\lambda_{ij})}{(\lambda_{ij}(q^2-1)-(q^2+1))}
      \frac{(q^2-1)(-1+\lambda_{ij})}{(\lambda_{ij}(q^2-1)-(q^2+1))}\]
\[ = \frac{2(q^4+1)-2\lambda_{ij} (q^4-1)  }{
            (\lambda_{ij}(q^2-1)-(q^2+1))^2}.   \]
Therefore $1-\beta_{ij} \gamma_{ij}=0$ if and only if
$\lambda_{i,j} = \frac{q^4+1}{q^4-1}$.

(ii) Using (11) we can replace every nonordered monomial by a sum
of lexicographically smaller monomials or monomials of the
type $v_j v_j v_i$.
Because of Lemma 2(i)
we can reduce $v_j v_j v_i$, $i<j$
to ordered monomials. Therefore we can reduce every
polynomial of degree three with finite steps to a
sum of ordered polynomials.
$\bullet$.

\begin{Le}
(i) Let $i<j<k$ be three indices,
$x_i,x_j,x_k$ three of the generators of $A$,
$X_{ij}\in I$ (cf. (7)),
$V\subset I$ the subspace of all elements
$x_a X_{ab}$, $x_b X_{ab}$, $X_{ab} x_a$, $X_{ab} x_b$ and
\[ X_{ijk}:=
-a_{jk} a_{ik} X_{ij}x_k-
 a_{jk}        x_jX_{ik}-
               X_{jk}x_i+
               x_kX_{ij}+
 a_{ij}        X_{ik}x_j+
 a_{ik} a_{ij} x_iX_{jk} \in I.   \]
Then we have
\[ X_{ijk}=
   \frac{(\lambda_{ij}\lambda_{jk}-\lambda_{ij}\lambda_{ik}
  -\lambda_{ik}\lambda_{jk} +1)(q^2-1)^2 (q^2+1)}{
                 k_{ij} k_{ik} k_{jk} } \times            \]
\begin{eqnarray}
\times
( p_1 x_i^3
+ p_2 x_j^3
+ p_3 x_k^3
+ p_4 x_i^2 x_j
+ p_5 x_i x_j^2+    \end{eqnarray}
\[p_6 x_i^2 x_k
+ p_7 x_i x_k^2
+ p_8 x_j^2 x_k
+ p_9 x_j x_k^2)\ \ mod(V)
\]
with
\[ p_1:=\frac{ 4(\lambda_{ij}-\lambda_{ik}) q^6 (q^2-1)}{
                             l_{ij} l_{ik}  }, \ \ \ \ \
   p_2:=\frac{ -2(q^2-1)K }{  l_{ij} l_{jk} }           \]
\[ p_3:=\frac{ 4(\lambda_{ik}-\lambda_{jk})      (q^2-1)}{
                             l_{ik} l_{jk}  }           \]
\[ p_4:=\frac{-2 m_{ij} q^2 }{  l_{ij}    },\ \ \ \
   p_5:=\frac{ 2 m_{ij}     }{  l_{ij}    },\ \ \ \
   p_6:=\frac{ 2 m_{ik} q^2 }{  l_{ik}    } \]
\[ p_7:=\frac{-2 m_{ik}     }{  l_{ik}    }\ \ \ \
   p_8:=\frac{-2 m_{jk} q^2 }{  l_{jk}    }\ \ \ \
   p_9:=\frac{ 2 m_{jk}     }{  l_{jk}    } \]
and
\[ k_{ij}:=1+q^2+\lambda_{ij} (1-q^2),\ \ \ \ \
   l_{ij}:=1+q^4+\lambda_{ij} (1-q^4),        \]
\[ m_{ij}:=1+q^2+\lambda_{ij} (-1+q^2),  \]
\[ K:=-(1+q^2+q^4+q^6)+\lambda_{12}( 1-q^2+q^4+q^6)
                      +\lambda_{23}(-1-q^2+q^4-q^6)+ \]
\[         \lambda_{12}\lambda_{23}( 1-q^2-q^4+q^6).  \]
The right side in $(12)$ vanishes if and only if
$\lambda_{ij}\lambda_{jk}-\lambda_{ij}\lambda_{ik}
   - \lambda_{ik}\lambda_{jk}= -1$  or
$\lambda_{ij}=\frac{1+q^2}{1-q^2}$, $\forall i<j$.
\end{Le}

\noindent
{\it Proof.}
One can check that in $X_{ijk}$ all terms $x_a x_b x_c$ with
three different $a,b,c$ cancel. The remaining monomials have
only two different indices.
Because we used in the reduction process of
$v_a v_b v_b$, $v_a v_a v_b$, $v_a v_b v_a$
in the proof of Lemma 2 only reduction rules which
correspond to
$x_a X_{ab}$, $x_b X_{ab}$, $X_{ab} x_a$, $X_{ab} x_b$,
we can replace
the remaining monomials of type $x_a x_b x_b$,
$x_a x_a x_b$, $x_a x_b x_a$ by sums of ordered polynomials
$mod(V)$. This was done, using a computer calculation.

The right side of (12) vanishes if and only if (2) is satisfied
or $p_1,...,p_9=0$.
From $p_1=p_3=p_4=0$ follows
$\lambda_{ij}=\lambda_{ik}=\lambda_{jk}=\frac{1+q^2}{1-q^2}$.
It follows $p_5,...,p_9=0$ and
$p_2=0$ because of
\[ (1-q^2)^2K= -(1-q^2)^2      (1+q^2+q^4+q^6)
               +(1-q^2)(1+q^2) (-2q^2+2q^4)
               +(1+q^2)^2      (1-q^2-q^4+q^6)=0.  \]
I.e. $p_1,...,p_9=0$ if and only if
$\lambda_{ij}=\frac{1+q^2}{1-q^2}$.
$\bullet$

\begin{Th}
Let ${A'}_\lambda$ be the algebra of $n$ points.
We can represent every polynomial of degree
less than three as a unique sum of ordered polynomials
$v_i,\ v_i v_j,\ v_i v_j v_k,\ \ i\leq j\leq k$ if and only if
the parameters
$\lambda_{ij}$, $\lambda_{ik}$, $\lambda_{jk}$
satisfy the conditions $(2)$:
\begin{eqnarray*}
     \lambda_{ij} \lambda_{jk} - \lambda_{ij} \lambda_{ik}
   - \lambda_{ik} \lambda_{jk} = -1, \ \ \ \ \  i < j < k
\end{eqnarray*}
and the additional conditions
$\lambda_{ij}, \lambda_{ik}, \lambda_{jk}
\neq \frac{q^2+1}{q^2-1}, \frac{q^4+1}{q^4-1},
\ \ \ \ \ i<j<k$ or if
\[\lambda_{ij}=\lambda_{ik}=\lambda_{jk}
=\frac{1+q^2}{1-q^2},\ \ \ \ \ i<j<k.\]
\end{Th}

\noindent
{\it Proof.}
(i)
Let $\lambda_{ij}=\frac{q^2+1}{q^2-1}$ (i.e. $d_{ij}= 0$).
Because of (7) we obtain the linear
dependence of ordered monomials
$a_{ij} v_i v_j + b_{ij} v_i^2 + c_{ij} v_j^2 = 0$.

Let $\lambda_{ij}=\frac{q^4+1}{q^4-1}$.
Because of Lemma 2 we obtain a linear
dependence of ordered monomials.

Finally let $i,j,k$ be three indices,
such that (2) is not satisfied and let
$\lambda_{ij}\neq\frac{1+q^2}{1-q^2}$ for $i<j$.
If we apply the quotient map $A \rightarrow A/I$ to (12)
we obtain a linear dependence between ordered monomials.

Therefore the conditions of Theorem 4 are necessary.

(ii)
Now let the conditions of the Theorem be satisfied.
Because the relations (9) are homogeneous,
${A'}_\lambda$ is the direct sum of the subspaces
of homogeneous polynomials.

Because of $\lambda_{ij}\neq \frac{q^2+1}{q^2-1}$
(i.e. $d_{ij}\neq 0$) we can represent every polynomial
of degree $\leq 2$ as a sum of ordered monomials.

It remains to consider monomials of degree three.
Because of $\lambda_{ij} \neq \frac{q^4+1}{q^4-1}$
and Lemma 2(ii)
we can reduce every monomial of degree three
in finite steps to a sum of ordered monomials.

We show that the ordered monomials are linear
independent.
The space of third degree polynomials in the
free algebra $A$ is $n^3$-dimensional.
The Ideal $I$ is generated by the
$2n\left( \begin{array}{c} n \\ 2 \end{array} \right)$
elements $x_i X_{jk}$, $X_{jk} x_i$, $j<k$.
Because (2) is satisfied or
$\lambda_{ij}=\frac{1+q^2}{1-q^2}$ for $i<j$,
the left side of (12) in Lemma 3 vanishes $mod(V)$.
Therefore we have
$\left( \begin{array}{c} n \\ 3 \end{array} \right)$
relations between the generators
$x_i X_{jk}$, $X_{jk} x_i$ of $I$.
These relations are independent, because
elements $x_i X_{jk}$, $X_{jk} x_i$
with three fixed different indices occur only in one
of the relations.
Therefore $dim\ I  \leq
 2n \left( \begin{array}{c} n \\ 2 \end{array} \right)-
    \left( \begin{array}{c} n \\ 3 \end{array} \right)$.
It follows, that the dimension of the space of third
degree polynomials in the factor algebra
${A'}_\lambda = A/I$
is greater than or equal
$n^3-2n\left( \begin{array}{c} n \\ 2 \end{array} \right)+
\left( \begin{array}{c} n \\ 3 \end{array} \right)
=n^2 +
\left( \begin{array}{c} n \\ 3 \end{array} \right)$.
We have shown above that this space is spanned
by ordered polynomials,
i.e. the dimension of this space is less than or equal
$n+2\left( \begin{array}{c} n \\ 2 \end{array} \right )
  +\left( \begin{array}{c} n \\ 3 \end{array} \right )
  =n^2 +
\left( \begin{array}{c} n \\ 3 \end{array} \right)$.
Therefore the ordered monomials are linear independent.
$\bullet$

\noindent
{\bf Remark:} If (2) is satisfied and
$\lambda_{ij} = \frac{q^2+1}{q^2-1}$
or $\frac{q^4+1}{q^4-1}$ for certain $i<j$
one can change the order of the indices in order to
satisfy the conditions of Theorem 4.

\noindent
{\bf Remark:}
We expect an analogue of Theorem 4 for polynomials of
higher degrees. One difficulty is to classify the
additional conditions.
For example, for fourth order polynomials
one can derive the additional condition
$\beta_{ij}\gamma_{ij}(1+\alpha_{ij})\neq 1$.

\noindent
{\bf Remark:}
If we introduce the new parameters
$c:= {\rm i} \frac{q^2 +1}{q^2-1}$ and
$\phi_{ij}:=\frac{c}{\lambda_{ij}}$,
($\lambda_{ij}\neq 0$)
we have
\[ [v_j,v_i]= \frac{{\rm i}}{c} (v_j^2-v_i^2)
             -\frac{{\rm i}}{\phi_{ij}} (v_i-v_j)^2. \]
For condition (2) we get the form
\[  \phi_{ik}= \frac{\phi_{ij} + \phi_{jk}}{
             1+\frac{\phi_{ij}   \phi_{jk}}{c^2} } \]
of the Einstein addition theorem for velocities.

In subsections C and D we give some examples.
\\ \\


\noindent
{\bf C. Points of the real projective line} \\
Consider $A_\lambda$ with real $\lambda_{ij}$
and endow $A_\lambda$ with the trivial involution $v_i^*=v_i$.
The relations (9) are  obviously $*$-invariant.
Therefore the involution on $A$ induces an involution
on ${A'}_\lambda$ and $A_\lambda$.

In view of Theorem 4 and the Remark after Theorem 4 we give the
following definition.

\begin{De} \rm
The  {\it algebra $A_\lambda^{\RR}$ of
coordinates of noncommutative points $v_i$
of the real line}
is the $U_q(sl(2,\RR ))$-modul $A_\lambda$,
if $A_\lambda$ carries the involution $*$ and
${A'}_\lambda$ has the property (2)
or is the exceptional structure
of Theorem 4.
\end{De}

\noindent
{\it Example $7$.} Let $I$ be an arbitrary linear ordered set
and let $\lambda_{ij}=1$ for $i<j$ (cf. Example 1).
It follows
\[  [v_j,v_i] = (q^2-1) v_i (v_j-v_i), \ \ \ i<j.  \]
We denote this structure by  $A^{\RR}_1$.

\noindent
{\it Example $8$.} Let $I$ be the set
$\ZZ$  of integers and
$ \lambda_{ij}= \coth ( \sum_{k=i}^{j-1} \alpha_k ),$ \\
$\alpha_k > 0,
\ \ \ i<j $ (cf. Example 2). It follows
\[ [v_j,v_i]=\frac{q^2-1}{q^2+1}
             (v_j^2-v_i^2 -
\coth ( \sum_{k=i}^{j-1} \alpha_k )(v_i-v_j)^2),
\ \ \ i<j. \]

\noindent
{\it Example $9$.} Let $I$ be the set of integers and
$\lambda_{ij}=\frac{1+q^2}{1-q^2}$, $i<j$
(cf. Theorem 4). It follows
\[ v_j v_i = \frac{ v_i^2 + q^2 v_j^2 }{ 1 + q^2 }.\]
We denote this structure by  $A^{\RR}_E$.
Because of $\lambda_{ij}\rightarrow \infty$
for $q\rightarrow 1$ this structure has no
quasiclassical limit.
\\ \\

\noindent
{\bf D. Points of the complexified projective line} \\
Let $I'$ be an arbitrary set and $I=I'\times \ZZ_2$.
We use the notations $w_i := v_{(i,0)}$ and
$ w^*_{{\overline{i}}} := v_{(i,1)}$, $i\in I'$.
Let $w_i,\ w_{\overline{i}}^*,\ i \in I'$ be connected
by the involution $*$
(i.e. $w_{\overline{i}}^{**}=w_i$).
Then the relations (9) get the form
\[ [w_j,w_i]=\frac{q^2-1}{q^2+1}
             (w_j^2-w_i^2 -\lambda_{ij}(w_i-w_j)^2),    \]
\[ [w_j,w_{\overline{i}}^*]=\frac{q^2-1}{q^2+1}
             (w_j^2-w_{\overline{i}}^{*2}-
           \lambda_{{\overline{i}}j}
             (w_{\overline{i}}^*-w_j)^2),    \]
\[ [w_{\overline{j}}^*,w_{\overline{i}}^*]=
    \frac{q^2-1}{q^2+1}
             (w_{\overline{j}}^{*2}-
              w_{\overline{i}}^{*2}-
              \lambda_{{\overline{i}}{\overline{j}}}
(w_{\overline{i}}^*-w_{\overline{j}}^*)^2).    \]
\begin{Le}
The above relations are $*$-invariant
(and therefore induce an involution on ${A'}_\lambda$),
if and only if
$\overline{\lambda_{ij}}
=\lambda_{ {\overline{i}} {\overline{j}} }$
and
$\overline{\lambda_{i {\overline{j}} }}
=\lambda_{{\overline{i}}j}$.
\end{Le}
We denote
the algebra generated by elements $w_i,w_i^*$ and
the above relations with
$\overline{\lambda_{ij}}
=\lambda_{ {\overline{i}} {\overline{j}} }$
and
$\overline{\lambda_{i {\overline{j}} }}
=\lambda_{{\overline{i}}j}$
by $A^c_\lambda$.

In order to check (2) for $A^c_\lambda$
it is sufficient to check
 $\lambda_{\alpha\beta} \lambda_{\beta\gamma}
-\lambda_{\alpha\beta} \lambda_{\alpha\gamma}
-\lambda_{\alpha\gamma}\lambda_{\beta\gamma}=-1$
for $(\alpha,\beta,\gamma)=
     (i,j,k),(\overline{i},j,k),(i,\overline{i},j)$,
$\forall i<j<k$.
(For real coefficients one has to check the case
$(\alpha,\beta,\gamma)= (i,j,k)$ and
$\lambda_{\overline{i}\overline{j}}=\lambda_{ij}$,
$\lambda_{i\overline{j}}=\lambda_{\overline{i}j}
=\lambda_{ij}^{-1}$ (i.e. $\lambda_{i\overline{i}}=0$).)

\begin{De} \rm  By the {\it algebra $A^{ \CC}_\lambda$ of
quantized complex points $w_i$ of the
complexified projective
line} we denote the $U_q(sl_2,\RR )$-modul $A^c_\lambda$,
if $A_\lambda^c$ satisfies (2) and
${A}^c_\lambda$ carries the involution $*$.
\end{De}
We consider the elements $w_i$ and  $w_{\overline{i}}^*$,
respectively,
as quantized complex  and
quantized complex conjugated points, respectively.

\noindent
{\it Example $10$.} Let $I'$ be a linear ordered set,
$I=I'\times \ZZ_2$
and let
$\lambda_{ij}$ be defined as in Example 4. We obtain

 \begin{eqnarray*}{}
  { [w_{\overline{i}}^*,w_{i}]}
  & = & \frac{q^2-1}{q^2+1}
                     (w_{\overline{i}}^{*2} - w_{i}^{2}),   \\
  { [w_j  ,w_i]  }   & = & (q^2-1)
            w_i (w_j-w_i),\ \ \ i<j,             \\
  { [w_{\overline{j}}^*,  w_{\overline{i}}^*] }  & = & (q^2-1)
     w_{\overline{i}}^*  (w_{\overline{j}}^{*}-w_{\overline{i}}^*),
    \ \ \ i<j,     \\
  { [w_j  ,w_{\overline{i}}^*]  } & = & (q^2-1)
           w_{\overline{i}}^*
          (w_j - w_{\overline{i}}^*),\ \ \ i<j,       \\
  { [w_{\overline{j}}^*,w_i  ]  } & = & (q^2-1)
            w_i (w_{\overline{j}}^{*}-w_i ),\ \ \ i<j.
\end{eqnarray*}
We denote this structure by $A^{ \CC}_1$

\noindent
{\it Example $11$.} Let $I'= \ZZ$,
$I=I'\times \ZZ_2$
and
$ \lambda_{ij}= \coth ( \sum_{k=i}^{j-1} \alpha_k ),\
\alpha_k > 0$, \\   $i<j $. It follows
\begin{eqnarray*}{}
{[w_{\overline{i}}^*,w_{i}]} & = & \frac{q^2-1}{q^2+1}
                    (w_{\overline{i}}^{2*} - w_{i}^{2}),   \\
 {[w_j,w_i]} & = & \frac{q^2-1}{q^2+1}
             (w_j^2-w_i^2 -
           \coth ( \sum_{k=i}^{j-1} \alpha_k )
(w_i-w_j)^2),\ \ \ i<j,                         \\
 {[w_j,w_{\overline{i}}^*]} & = & \frac{q^2-1}{q^2+1}
             (w_j^2-w_{\overline{i}}^{*2}-
           \tanh ( \sum_{k=i}^{j-1} \alpha_k )
             (w_{\overline{i}}^*-w_j)^2),\ \ \ i<j,         \\
 {[w_{\overline{j}}^*,w_{i}]} & = & \frac{q^2-1}{q^2+1}
             (w_{\overline{j}}^{*2}-w_{i}^{2}-
           \tanh ( \sum_{k=i}^{j-1} \alpha_k )
             (w_{i}-w_{\overline{j}}^*)^2),\ \ \ i<j,         \\
 {[w_{\overline{j}}^*,w_{\overline{i}}^*]}
 & = & \frac{q^2-1}{q^2+1}
             (w_{\overline{j}}^{*2}-
              w_{\overline{i}}^{*2}-
           \coth ( \sum_{k=i}^{j-1} \alpha_k )
(w_{\overline{i}}^*-w_{\overline{j}}^*)^2),\ \ \ i<j.
\end{eqnarray*}
Example 10 arises again in the limit $\alpha_k \rightarrow \infty$.

\noindent
{\it Example $12$.}
An example with complex coefficients
is given by
$I=\{1,2,\overline{1},\overline{2} \}$ and
$\lambda_{12}=\lambda_{2\overline{2}}
=\lambda_{\overline{2}\overline{1}}
=\lambda_{\overline{1}1}={\rm i}$
and $\lambda_{1\overline{2}}
    =\lambda_{\overline{1}2}=0$
(cf. Example 6).
\\ \\

\noindent
{\bf E. Algebras of three points}

\noindent
We will say that ${A_\lambda}'$ is polynomial, if
the subspaces of homogeneous polynomials have the classical
dimensions.

From Theorem 4 and the Remark after Theorem 4 follows
that the subspaces of homogeneous polynomials of degree $\leq 3$
of the algebras ${A_\lambda}'$ which satisfy (2)
and of the exceptional algebra ${A_E^\RR}$
have the classical dimensions.

We expect that these algebras are polynomial.
We give a proof for algebras of three points.

\begin{Le}
Let ${A_\lambda}'$ be the algebra with three generators
$v_1,v_2,v_3$ and let $(2)$ be satisfied.
${A_\lambda}'$ is equivalent to the algebra with generators
$u_1,u_2,u_3$ and relations
\[   u_2 u_1 =  \frac{1}{q^2}  u_1 u_2,  \]
\begin{eqnarray}
     u_3 u_1 =  \frac{1}{q^2} u_1 u_3
+(\frac{1}{q^2}-1)(\lambda_{1,2}-1)(\lambda_{2,3}+1) u_2^2,
\end{eqnarray}
\[   u_3 u_2 =  \frac{1}{q^2} u_2 u_3  \]
The equivalence is given by a linear transformation
between the generators $u_i$ and $v_i$.
\end{Le}

\noindent
{\it Proof.}
Consider the transformation
{\tiny
\begin{eqnarray}
\left( \begin{array}{c}
u_1 \\
u_2 \\
u_3
\end{array} \right)
=
\left( \begin{array}{ccc}
\lambda_{1,2} -1    & -(\lambda_{1,2}+\lambda_{2,3}) & \lambda_{2,3} +1 \\
-(\lambda_{1,2} +1) & \lambda_{1,2}+\lambda_{2,3} & -(\lambda_{2,3} -1) \\
  (\lambda_{1,2}+1)^2 (\lambda_{2,3}+1) &
 -(\lambda_{1,2}+\lambda_{2,3})(\lambda_{1,2}-1)(\lambda_{2,3}+1) &
  (\lambda_{2,3}-1)^2 (\lambda_{1,2}-1)
\end{array} \right)
\left( \begin{array}{c}
v_1 \\
v_2 \\
v_3
\end{array} \right).
\end{eqnarray}
}
We denote the matrix by $(k_{ij})$.
For the determinant one obtains
\[ det((k_{ij})) =- 8 (\lambda_{1,2}+\lambda_{2,3})^2.  \]
(If $det((k_{ij}))= 0$
there exists a permutation of the three indices such that
$det((k_{ij}))\neq 0$.)

Let
\[   A_{12}:= 2 u_1 u_2 - 2 q^2 u_2 u_1  ,  \]
\begin{eqnarray}
     A_{13}:= 2 u_1 u_3 - 2 q^2 u_3 u_1
+(1-q^2)(\lambda_{1,2}-1)(\lambda_{2,3}+1) u_2^2 ,
\end{eqnarray}
\[   A_{23}:= 2 u_2 u_3 - 2 q^2 u_3 u_2  .  \]
Because of (13) we have $A_{12}=A_{13}=A_{23}=0$.
We form the vanishing expressions
\[  B_{12}: =
 ( k_{33} A_{12} -k_{23} A_{13} +k_{13} A_{23} )/det(K), \]
\begin{eqnarray}
    B_{13}: =
 (-k_{32} A_{12} +k_{22} A_{13} -k_{12} A_{23} )/det(K),
\end{eqnarray}
\[  B_{23}: =
 ( k_{31} A_{12} -k_{21} A_{13} +k_{11} A_{23} )/det(K). \]
We insert (14). We obtain with a computer calculation
\[  B_{12}=
 -(1+q^2) (v_2 v_1- v_1 v_2) +
 (1-q^2) (v_1^2-v_2^2+\lambda_{12}(v_1-v_2)^2)=0 \]
\[  B_{13}=
 -(1+q^2) (v_3 v_1- v_1 v_3) +
 (1-q^2)
(v_1^2-v_3^2+
\frac{1+\lambda_{12}\lambda_{23}}{\lambda_{12}+\lambda_{23}}
(v_1-v_3)^2) =0     \]
\[ B_{23}=
 -(1+q^2) (v_3 v_2- v_2 v_3) +
 (1-q^2) (v_2^2-v_3^2+\lambda_{13}(v_2-v_3)^2) =0 \]
Using (2) the formulas (9) follow.
$\bullet$

\begin{Sa} Let $(2)$ be satisfied.
Then ${A_\lambda}'$ is polynomial.
\end{Sa}

\noindent
{\it Proof.}
The transformation between the $v_i$ and $u_i$ does not
change the dimensions
of the subspaces of homogeneous polynomials.
We have to show the polynomiality of the algebra (13).

By the Diamond lemma (cf. Ref. 4)
this algebra has the PBW-property
(and is therefore "polynomial")
if and only if the "overlap"
\[ (u_3 u_2) u_1 - u_3 (u_2 u_1) \]
gives zero when it is reduced by means of (13) to a linear
combination of ordered monomials (cf. Ref. 4 for details).
Let $a:=\frac{1}{q^2}$ and
$\epsilon:=(\frac{1}{q^2}-1)(\lambda_{1,2}-1)(\lambda_{2,3}+1)$.
We obtain
\[   (u_3 u_2) u_1 - u_3 (u_2 u_1)            \]
\[ = a u_2 (u_3 u_1) - a (u_3 u_1) u_2        \]
\[ =   a^2 (u_2 u_1) u_3 + a \epsilon u_2^3
     - a^2 u_1 (u_3 u_2) - a \epsilon u_2^3   \]
\[ = a^3 u_1 u_2 u_3 - a^3 u_1 u_2 u_3 = 0.\bullet \]
\\ \\

Finally we consider the exceptional structure $A^{\RR}_E$.
It is given by
\begin{eqnarray}
v_j v_i = \frac{ v_i^2 + q^2 v_j^2 }{ 1 + q^2 },\ \ \ \ \
i<j,
\end{eqnarray}
(cf. Example 9).
\begin{Sa}
The Algebra $A_E^\RR$ is polynomial.
\end{Sa}

\noindent
{\it Proof.}
Consider the transformation
\begin{eqnarray}
\left( \begin{array}{c}
v_1 \\
v_2 \\
v_3
\end{array} \right)
=
\left( \begin{array}{ccc}
q^2 & q^2 & q^2 \\
1   & q^2 & q^2 \\
1   & 1   & q^2
\end{array} \right)
\left( \begin{array}{c}
u_1 \\
u_2 \\
u_3
\end{array} \right).
\end{eqnarray}
We denote the matrix by $(k_{ij})$.
We obtain for the determinant
\[ det((k_{ij})) = q^2(q^2-1)^2 \neq 0.  \]
Let
\[   A_{12}:=-(1+q^2)v_2 v_1+v_1^2+q^2 v_2^2  ,  \]
\begin{eqnarray}
     A_{13}:=-(1+q^2)v_3 v_1+v_1^2+q^2 v_3^2  ,
\end{eqnarray}
\[   A_{23}:=-(1+q^2)v_3 v_2+v_2^2+q^2 v_3^2  .  \]
Because of (17) we have $A_{12}=A_{13}=A_{23}=0$.
We form the vanishing expressions $B_{12}$, $B_{13}$, $B_{23}$
(cf. (16)) and insert (19). A calculation yields
\[  B_{12}= -(1+q^2) u_2 u_1 = 0,                      \]
\[  B_{13}= u_1 u_2 + u_2 u_1 + u_1 u_3 - q^2 u_3 u_1 = 0, \]
\[  B_{23}=-u_1 u_2 - u_2 u_1 + u_2 u_3 - q^2 u_3 u_2 = 0. \]
Therefore the algebra has generators $u_1,u_2,u_3$ and relations
\[  u_2 u_1:=0, \]
\begin{eqnarray}
    u_3 u_1:=\frac{1}{q^2} (u_1 u_3 + u_1 u_2),
\end{eqnarray}
\[  u_3 u_2:=\frac{1}{q^2} (u_2 u_3 - u_1 u_2). \]
The monomials on the left side are lexicographically greater
than the monomials on the right side.
Therefore it is sufficient to show, that the "overlap"
\[ (u_3 u_2) u_1 - u_3 (u_2 u_1) \]
gives zero when it is reduced by means of (20) to a linear
combination of ordered monomials (cf. Proof of Proposition 7).
We have
\[ (u_3 u_2) u_1 - u_3 (u_2 u_1) \]
\[ =\frac{1}{q^2} (u_2 (u_3 u_1) - u_1 (u_2 u_1) ) \]
\[ =\frac{1}{q^4} ((u_2 u_1) u_3 + (u_2 u_1) u_2) = 0.\bullet \]
\\

\noindent
{\bf F. Projective Coordinates} \\
As the algebra of {\it noncommutative projective
coordinates $B_\mu$
for $A_\lambda$} we denote the
$U_q(sl(2,\RR ))$-modul algebra generated
by elements $x_i,y_i,y_i^{-1},
(\sum \alpha_i x_i y_i^{-1} )^{-1}$, $i,j\in I$
and relations

\begin{eqnarray*}{}
            y_i y_i^{-1}   & = & y_i^{-1}y_i             = 1, \\
      (\sum \alpha_i x_iy_i^{-1})^{-1}
      (\sum \alpha_i x_iy_i^{-1})      &
  = & (\sum \alpha_i x_iy_i^{-1})
      (\sum \alpha_i x_iy_i^{-1})^{-1}  = 1
   \end{eqnarray*}
and
 \begin{eqnarray}
  x_iy_i             &  = &  qy_ix_i,                 \\
  x_jx_i             &  = & \mu_{ij}^{(1)} x_i x_j,      \\
  y_jy_i             &  = & \mu_{ij}^{(1)} y_i y_j,      \\
  x_jy_i             &  = & \mu_{ij}^{(2)} y_i x_j
   + (\mu_{ij}^{(1)} - \frac{1}{q} \mu_{ij}^{(2)} ) x_i y_j, \\
  y_jx_i             &  = & \mu_{ij}^{(2)} x_i y_j
   + (\mu_{ij}^{(1)} -    q \mu_{ij}^{(2)}  )       y_i x_j
\end{eqnarray}
with
\[ \mu_{ji}^{(1)} = \frac{ 1 }{ \mu_{ij}^{(1)} }
   \ \ \ \ \    {\rm and}  \ \ \ \ \
   \mu_{ji}^{(2)} = \frac{ 1 }{ \mu_{ij}^{(2)}
            (\frac{ \mu_{ij}^{(1)} }{ \mu_{ij}^{(2)} }
            (q+\frac{1}{q}
           -\frac{ \mu_{ij}^{(1)} }{ \mu_{ij}^{(2)} }))}. \]
The $U_q(sl(2,\RR ))$-modul structure on $B_\mu$ is given by
\begin{eqnarray*}
 E(x_i)    =  q^\frac{1}{2} y_i,             &
 \ E(y_i)  =  0,                          \\
 F(x_i)    =  0,                          &
 \ F(y_i)  =  q^{-\frac{1}{2}}x_i,        \\
 K(x_i)    =  q^{-\frac{1}{2}}x_i ,          &
 \ K(y_i)  =  q^\frac{1}{2}   y_i.
\end{eqnarray*}
and equations (4), (5), (6).

\begin{Le}
The elements
$v_i: = q^{  \frac{1}{2} } x_i y_i^{-1}$, $i\in I$ satisfy $(9)$
with
\[ \lambda_{ij} =
   \frac{q^2-2q \frac{ \mu^{(1)}_{ij}}{ \mu^{(2)}_{ij} }
         +1}{1-q^2} \ \ \ \ \
{\rm or} \ \ \ \ \
\mu_{ij}:=\frac{\mu^{(1)}_{ij}}{\mu^{(2)}_{ij}}=
 \frac{ 1 + q^2 + \lambda_{ij} (q^2-1) }{ 2q }.\]
and $U_q(sl(2,\RR))$ acts on $v_i$ according to $(3)$.
\end{Le}
This can be checked by an explicite calculation.
Therefore we can map $A_\lambda$ onto a submodul of
$B_\mu$.
We call $x_i,y_i$ {\it noncommutative coordinates} of $v_i$.

\noindent
{\bf Remark:} One can show
that the algebra of
projective coordinates is {\it polynomial}
if and only if
$\mu^{(1)}_{ij}=q^{\pm 1}\mu^{(2)}_{ij}$, i.e.
$\lambda_{ij}=\pm 1$.
Therefore for $\lambda_{ij}\neq \pm 1$ the
polynomiality of the algebra  ${A'}_\lambda$ does not yield
the polynomiality of the algebra of projective coordinates.

\newpage
\noindent
{\Large \bf IV. Quantum cross ratios }
\\ \\
\noindent
{\bf A. Invariant elements} \\
In this section we investigate
a class of invariant elements
in $A_\lambda$.

\begin{De}\rm We say
that $x \in A_\lambda$ (or $B_\mu$)
is {\it invariant}, if
\[ E(x)=0,\ \ \ \ \  F(x)=0,\ \ \ \ \  K(x)=x. \]
\end{De}

\noindent
Following Ref. 1 we define
\[ (ij):= q^{-\frac{1}{2}}  x_i y_j -
q^{\frac{1}{2}} y_i x_j ,\ \ \ \ \
{\rm and} \ \ \ \ \ [ij]:=v_i-v_j. \]
One can check that $(ij)$ is invariant and we have
\begin{eqnarray}
      (ij) =  {y_i} [ij] {y_j}.
\end{eqnarray}
Now we use the invariants $(ij)\in B_\mu$ to combine
invariants of $A_\lambda$.
We define
\[  C_{ijkl}:= (il)(kl)^{-1} (kj) (ij)^{-1}. \]

\begin{Sa} (i) $C_{ijkl}\in B_\mu $, ($i,j,k,l\in I$)
is invariant. \\
(ii) We can represent $C_{ijkl}$ by elements $v_i$
(i.e. $C_{ijkl}\in A_\lambda$)
and we have the formula
{\small \[  C_{ijkl}:=
       (q+{q}^{-1}  - \mu_{il}) (v_i-v_l)
     ( (q+{q}^{-1}  - \mu_{ik}) v_k
     - (q+{q}^{-1}  - \mu_{il}) v_l
     + (\mu_{ik} - \mu_{il}) v_i  )^{-1} \times \]
\begin{eqnarray}
     ((q+{q}^{-1}  - \mu_{ik}) v_k
     -(q+{q}^{-1}  - \mu_{ij}) v_j
     +( \mu_{ik} -\mu_{ij})v_i  )
      (q+{q}^{-1}  -\mu_{ij})^{-1} (v_i-v_j)^{-1}.
\end{eqnarray}    }
\end{Sa}
{\it Proof.}
The element $(ij)^{-1}$ and
products of invariant elements are invariant
because of (4), (5), (6). The assertion follows.

\noindent
(ii) From (21)-(25) follow the commutation rules
\begin{eqnarray}
y_i v_i = q^{-1}  v_i y_i, \ \ \ \ \
y_i v_j = ((q+q^{-1}-\mu_{ij}) v_j +(\mu_{ij}-q) v_i)y_i.
\end{eqnarray}
Further, because of (26) we have
\begin{eqnarray}
C_{ijkl}=y_i [il][kl]^{-1} [kj] [ij]^{-1} y_i^{-1}.
\end{eqnarray}
The assertion follows, if we move $y_i$ to the right
applying (28). $\bullet$

\begin{Sa} In analogy to the classical case we have
\[ C_{ilkj} = \frac{1}{C_{ijkl}}, \ \ \ \ \
   C_{ijlk} = 1-C_{ijkl},         \ \ \ \ \
   C_{ikjl} = \frac{C_{ijkl}}{C_{ijkl}-1},   \]
\[ C_{iklj} = \frac{1}{1-C_{ijkl}}, \ \ \ \ \
   C_{iljk} = 1-\frac{1}{C_{ijkl}},          \]
\end{Sa}

\noindent
{\it Proof.} One checks the equations
$C_{ilkj} C_{ijkl} = 1$ and
$C_{ijlk}+C_{ijkl} = 1$ using (29) and the rule
$[ik]=[ij]+[jk]$.
The five identities
can  combined from them. $\bullet$

\noindent
{\bf Remark:}
In the general case, the elements $C_{ijkl}$
are not $*$-invariant and the identities
$C_{ijkl}=C_{lkji}=C_{klij}=C_{jilk}$
are in general  not valid (see below).
\\ \\

\noindent
{\bf B. The cross ratio} \\
\noindent
We restrict the consideration to the structure $A_1^{\RR}$ of Example 7,
i.e. $\lambda_{ij}=1$.
(If $|\lambda_{ij}|=1$, we can reorder the indices such that
$\lambda_{ij}=1$ (cf. Example 1).)

We can choose the algebra of projective coordinates
with $\mu^{(1)}_{ij}=q^{2}$,
     $\mu^{(2)}_{ij}=q$      (cf. Lemma 9), i.e.
\[ x_i y_i               =   q      y_i x_i, \]
\[ x_j x_i               =   q^{2}  x_i x_j, \]
\[ y_j y_i               =   q^{2}  y_i y_j, \]
\[ x_j y_i               =   q      y_i x_j
                         +(q^{2}-1) x_i y_j,  \]
\begin{eqnarray}
   y_j x_i               =   q      x_i y_j.
\end{eqnarray}

\begin{Le}
Let $i<j<k<l$  be four ordered indices. Then we have
\[  (kl)(ij)  =   q^{ 6} (ij)(kl),   \]
\[  (jk)(il)  =          (il)(jk),   \]
\[  (jl)(ik)  =   q^{ 4} (ik)(jl) + (q^{ 6}-q^{ 4}) (ij)(kl),\]
\[  (jk)(ij)  =   q^{ 4} (ij)(jk),   \]
\[  (ik)(ij)  =   q^{ 2} (ij)(ik),   \]
\[  (jk)(ik)  =   q^{ 2} (ik)(jk),   \]
\[  (ij)^*    =           (ij).       \]
\end{Le}
\noindent
{\it Proof.} We prove the fourth rule.
Let $i<j<k$. From (30) follows
\[        x_k (x_i y_j - q y_i x_j)
   = q^2 x_i x_k y_j - q^2 y_i x_k x_j - (q^3-q) x_i y_k x_j \]
\[ = q^3 x_i y_j x_k + (q^4-q^2)x_i x_j y_k
   - q^4 y_i x_j x_k - (q^4-q^2)x_i x_j y_k \]
\[=  q^3 (x_i y_j - q y_i x_j) x_k,  \]
i.e. $x_k (ij)= q^3 (ij) x_k$.
Similar one derives
$y_k (ij)= q^3 (ij) y_k$,
$x_j (ij)= q   (ij) x_j$,
$y_j (ij)= q   (ij) y_j$.
The identity $(jk)(ij) = q^4 (ij)(jk)$ follows.
The proof of the remaining rules is similar.
$\bullet$

In a similar situation
some of the invariant
elements $C_{ijkl}$
were first considered in Ref. 1.
Let $i<j<k<l$ be four elements of $I$
in increasing order.
With $C:=C_{ijkl}$ we derive using Proposition 11 and Lemma 12
\newpage
\begin{eqnarray*}
C_{ijkl}=C_{klij}= C,        &       \ \ \ \ \
C_{jilk}=C_{lkji}= C^* ={q^2}C.
\end{eqnarray*}
\begin{eqnarray*}
C_{ilkj}=C_{kjil}= \frac{1}{C},&    \ \ \ \ \
C_{jkli}=C_{lijk}= C_{ilkj}^* = \frac{1}{q^2 C}
= \frac{ C_{ilkj} }{ q^2 },
\end{eqnarray*}
\begin{eqnarray*}
C_{ijlk}=C_{klji}= 1-C,         &   \ \ \ \ \
C_{jikl}=C_{lkij}= C_{ijlk}^* = 1-{q^2}C
= {q^2}C_{ijlk}+(1-{q^2}),
\end{eqnarray*}
\begin{eqnarray*}
C_{ikjl}=C_{kilj}= \frac{C}{C-1},&  \ \ \ \ \
C_{jlik}=C_{ljki}= C_{ikjl}^* = \frac{ q^2 C }{ q^2 C - 1 }
= \frac{ q^2 C_{ikjl} }{ 1 + (q^2 -1 ) C_{ikjl} },
\end{eqnarray*}
\begin{eqnarray*}
C_{iklj}=C_{kijl}= \frac{1}{1-C},&  \ \ \ \ \
C_{jlki}=C_{ljik}= C_{iklj}^* = \frac{ 1 }{ 1 - q^2 C }
= \frac{  C_{iklj} }{ q^2 + ( 1 - q^2 ) C_{iklj} },
\end{eqnarray*}
\begin{eqnarray*}
C_{iljk}=C_{kjli}=1-\frac{1}{C},&  \ \ \ \ \
C_{jkil}=C_{likj}=C_{iljk}^*= 1-\frac{1}{ q^2 C }
= \frac{1}{ q^2 } C_{iljk}+(1- \frac{1}{ q^2 } ).
\end{eqnarray*}
The Remark after Proposition 11
indicates, that $C_{ijkl}$ is not a good candidate
of a Quantum cross ratio.
Now we make the following observation:

\begin{Le}
Let $i,j,k,l$ be four indices in arbitrary order
and let $C_{ijkl}^*=:f(C_{ijkl},q^2)$,
where $f$ is a certain rational function.
Then $f( C_{ijkl} , q )$
is a $*$-invariant element of $A^{\RR}_\lambda$.
\end{Le}
This can be checked with the above formulas.
Now we give the following definition.
\begin{De}\rm We say that the $*$-invariant
element $f( C_{ijkl} , q )$ is the
{\it Quantum cross ratio} of
the Quadrupel $(v_i, v_j, v_k, v_l)$.
We use the notation  $ {\cal C}_{ijkl} $.
\end{De}
The elements ${\cal C}_{ijkl}$
have the desired properties. They tend to the classical
cross ratio for $q\rightarrow 1$ and for
$i<j<k<l$ we obtain
\newpage
\begin{eqnarray*}
{\cal C}_{ijkl} = {qC}{} =q [il][kl]^{-1}[kj][ij]^{-1} ,
 \end{eqnarray*}
 \begin{eqnarray*}
{\cal C}_{ilkj} = \frac{1}{ {\cal C}_{ijkl}  }
                   = \frac{1}{   q    C  }
                   = \frac{1}{q}      C_{ilkj},
 \end{eqnarray*}
 \begin{eqnarray*}
{\cal C}_{ijlk} = 1-  { \cal C }_{ijkl} =
                     1-{ q C }
                   = {q} C_{ ijlk } + ( 1- {q} ),
 \end{eqnarray*}
 \begin{eqnarray*}
{\cal C}_{ikjl} = \frac{{\cal C}_{ijkl}}{{\cal C}_{ijkl}-1} =
                     \frac{ q C }{ q C - 1 }
                   = \frac{ q C_{ikjl} }{ 1 + ( q-1 ) C_{ikjl} },
 \end{eqnarray*}
 \begin{eqnarray*}
{\cal C}_{iklj} = \frac{1}{1-{\cal C}_{ijkl}} =
                     \frac{1}{ 1- q C }
                   = \frac{ C_{iklj} }{ q + ( 1-q  ) C_{iklj}},
 \end{eqnarray*}
 \begin{eqnarray*}
{\cal C}_{iljk} = 1-\frac{1}{ {\cal C}_{ijkl} }
                   = 1-\frac{ 1 }{ q C }
                   = \frac{1}{q} C_{iljk}+( 1 - \frac{1}{q} ),
\end{eqnarray*}
and
\[ {\cal C}_{ijkl} = {\cal C}_{lkji} =
   {\cal C}_{klij} = {\cal C}_{jilk}.   \]
\\ \\

\noindent
{\bf C. The Euclidean distance} \\
In the classical case an Euclidean distance on $P\RR$
can be defined, if we fix the three points $t_0$ (zero-point),
$t_1$ (one-point) and $t_\infty$ (the point at infinity).
The distance $d(t_j,t_i)$ of two points $t_j,t_i$ is given by
\[ d(t_j,t_i) := C_{0,1,\infty,j} - C_{0,1,\infty,i}. \]
This observation leads us to a notion of a noncommutative
Euclidean distance. Fix three different generators
$v_0$, $v_1$ and $v_\infty$ of $A_1^{\RR}$. Then the
{\it noncommutative distance} of the
generators $v_i$ and $v_j$
is given by the $*$-invariant element
\[ D_{v_0,v_1,v_\infty} (v_j,v_i):=
 {\cal C}_{0,1,\infty,j} - {\cal C}_{0,1,\infty,i}. \]

 \noindent
 \setcounter{equation}{0}
 \vspace{0.5cm}

\small

\noindent
${}^1 \rm C$.-S. Chu, P.-M. Ho, B. Zumino,
"The Braided Quantum 2-Sphere", \\
$\ \ \  $ preprint q-alg 9507013,
UCB-PTH-95/25.

\noindent
${}^2 \rm J$.-H. Lu, A. Weinstein,
J. Diff. Geom. {\bf 31}, 501 (1990).

\noindent
${{ }^3 \rm S.}$ Klimek, A. Lesniewski,
{ J. Func. Anal.} {\bf 115}, 1 (1993).

\noindent
${{ }^4 \rm G.}$M. Bergman,
Adv. in Math. {\bf 29} (1978), 178-218.

\end{document}